\documentclass[a4paper,12pt]{article}
\usepackage{amsmath,amssymb,amsfonts}
\usepackage[all]{xy}
\def\mar[#1]{\ar@{-}[#1]|-{\object@{<}}}

\newcommand\boite[1]{\bullet\makebox[0pt][l]{$\scriptstyle #1$}}
\newcommand\boiteb[1]{\circ\makebox[0pt][l]{$\scriptstyle #1$}}
\CompileMatrices

\newcommand{\gmod}[1]{#1\hbox{-}\mathsf{Mod}}
\newcommand{\ec}[1]{\End_{\mathcal{C}}(#1)}
\newcommand{\hcp}[2]{\Hom_{\mathcal{C}_{p}}(#1,#2)}
\newcommand{\hc}[2]{\Hom_{\mathcal{C}}(#1,#2)}

\newcommand{\estliemod}[1]{\;\raisebox{.5ex}{\rule{2.5ex}{.2ex}}_{#1}\;}
\def\un{{\bf 1}}
\def\zero{\{0\}}
\def\mpoint{\;\;\;.}
\def\mvirg{\;\;\;,}
\def\mpn{\medskip\par\noindent}
\def\smp{\smallskip\par}
\def\normal{\mathop{\underline\triangleleft}}
\newcommand{\scal}[2]{\langle #1,#2\rangle}

\def\Id{\hbox{\rm Id}}

\def\Res{\hbox{\rm Res}}
\def\Ind{\hbox{\rm Ind}}
\def\Hom{\hbox{\rm Hom}}
\def\End{\hbox{\rm End}}

\def\Inf{\hbox{\rm Inf}}
\def\Def{\hbox{\rm Def}}
\def\Ten{\hbox{\rm Ten}}
\def\Iso{\hbox{\rm Iso}}
\def\Im{\hbox{\rm Im}}
\def\Ker{\hbox{\rm Ker}}

\def\Defres{\hbox{\rm Defres}}
\def\Teninf{\hbox{\rm Teninf}}
\def\Indinf{\hbox{\rm Indinf}}

\def\op{^{op}}

\def\dom{\backslash}

\newcommand{\sur}[1]{\,\overline{\! #1}}
\newcommand{\sumb}[2]{\sum_{{\scriptstyle #1}\atop {\scriptstyle #2}}}

\newcommand{\ressort}[1]{\hskip #1 plus #1 minus #1}
\def\findemo{~\leaders\hbox to 1em{\hss\ \hss}\hfill~\raisebox{.5ex}{\framebox[1ex]{}}\smp}

\makeatletter
\renewenvironment{enumerate}{\ifnum \@enumdepth >3 \@toodeep\else
      \advance\@enumdepth \@ne
      \edef\@enumctr{enum\romannumeral\the\@enumdepth}\list
      {\csname label\@enumctr\endcsname}{\setlength{\topsep}{1ex}\setlength{\itemsep}{0pt}\usecounter
        {\@enumctr}\def\makelabel##1{\hss\llap{##1}}}\fi}{\endlist}
\renewenvironment{itemize}{\ifnum \@itemdepth >3 \@toodeep\else \advance\@itemdepth \@ne
\edef\@itemitem{labelitem\romannumeral\the\@itemdepth}%
\list{\csname\@itemitem\endcsname}{\setlength{\topsep}{1ex}\setlength{\itemsep}{0pt}\def\makelabel##1{\hss\llap{##1}}}\fi}
{\endlist}

\def\@sect#1#2#3#4#5#6[#7]#8{\ifnum #2>\c@secnumdepth
     \let\@svsec\@empty\else
     \refstepcounter{#1}\edef\@svsec{\csname the#1\endcsname .\hskip .5em}\fi
     \@tempskipa #5\relax
      \ifdim \@tempskipa>\z@
        \begingroup #6\relax
          \@hangfrom{\hskip #3\relax\@svsec}{\interlinepenalty \@M #8\par}%
        \endgroup
       \csname #1mark\endcsname{#7}\addcontentsline
         {toc}{#1}{\ifnum #2>\c@secnumdepth \else
                      \protect\numberline{\csname the#1\endcsname}\fi
                    #7}\else
        \def\@svsechd{#6\hskip #3\relax  
                   \@svsec #8\csname #1mark\endcsname
                      {#7}\addcontentsline
                           {toc}{#1}{\ifnum #2>\c@secnumdepth \else
                             \protect\numberline{\csname the#1\endcsname}\fi
                       #7}}\fi
     \@xsect{#5}}

\def\section{\pagebreak[3]\setcounter{prop}{0}\setcounter{equation}{0}\@startsection{section}{1}{\z@}{4ex plus 4ex}{4ex}{\center\reset@font\large\bf}}
\def\subsection{\pagebreak[3]\refstepcounter{prop}\@startsection{subsection}{2}{\z@}{4ex plus 6ex}{-1em}{\reset@font\bf}}
\def\subsubsection{\@startsection{subsubsection}{3}{\z@}{4ex plus 6ex}{-1em}{\reset@font\it}}
\makeatother

\def\Z{\mathbb{Z}}
\def\1{{1\;\!\!\!{\rm l}}}

\def\Q{\mathbb{Q}}

\def\F{\mathbb{F}}

\def\theprop{\thesection.\arabic{prop}}
\renewenvironment{equation}{\refstepcounter{prop}$$}{\leqno{(\thesection.\arabic{prop}})$$}

\renewenvironment{equation}{\refstepcounter{subsection}\refstepcounter{prop}$$}{\leqno{\bf (\theprop)}$$}

\def\pf{\noindent{\bf Proof: }}
\newenvironment{rem}[1]{\refstepcounter{subsection}\refstepcounter{prop}\mpn{{\bf \thesection.\arabic{prop}.}\ \ \bf#1~:}}{\smp}
\newenvironment{enonce}[1]{\pagebreak[3]\refstepcounter{subsection}\refstepcounter{prop}\mpn{{\bf \thesection.\arabic{prop}.\ \ #1~:}}\begin{it} }{\end{it}\smp}

\def\thesection{\arabic{section}}

\begin{document}
\centerline{\Large\bf The functor of units of Burnside rings for $p$-groups}
\vspace{1cm}
\centerline{\bf Serge Bouc}
\vspace{.5cm}
\centerline{\footnotesize LAMFA - UMR 6140 - Universit\'e de Picardie-Jules Verne}
\centerline{\footnotesize 33 rue St Leu - 80039 - Amiens Cedex 1 - France}
\centerline{\footnotesize\tt email : serge.bouc@u-picardie.fr}
\def\thefootnote{}\def\thefootnote{\arabic{footnote}}
\vspace{1cm}
{\footnotesize\bf Abstract:} {\footnotesize In this note I describe the structure of the biset functor $B^\times$ sending a $p$-group $P$ to the group of units  of its Burnside ring $B(P)$. In particular, I show that $B^\times$ is a rational biset functor. It follows that if $P$ is a $p$-group, the structure of $B^\times(P)$ can be read from a genetic basis of $P$~: the group $B^\times(P)$ is an elementary abelian 2-group of rank equal to the number isomorphism classes of rational irreducible representations of~$P$ whose type is trivial, cyclic of order 2, or dihedral.}
\section{Introduction}
If $G$ is a finite group, denote by $B(G)$ the Burnside ring of $G$, i.e. the Grothendieck ring of the category of finite $G$-sets (see e.g. \cite{handbook}). The question of structure of the multiplicative group $B^\times(G)$ has been studied by T.~tom~Dieck~(\cite{tomdieckgroups}), T.~Matsuda~(\cite{matsuda}), T.~Matsuda and T.~Miyata~(\cite{matsudamiyata}), T.~Yoshida~(\cite{yoshidaunit}), by geometric and algebraic methods.\par
Recently, E. Yal\c cin wrote a very nice paper (\cite{yalcin}), in which he proves an induction theorem for $B^\times$ for $2$-groups, which says that if $P$ is a $2$-group, then any element of $B^\times(P)$ is a sum of elements obtained by inflation and tensor induction from  sections $(T,S)$ of $P$, such that $T/S$ is trivial or dihedral.\par
The main theorem of the present paper implies a more precise form of Yal\c cin's Theorem, but the proof is independent, and uses entirely different methods. In particular, the biset functor techniques developed in \cite{doublact}, \cite{fonctrq} and~\cite{bisetsections}, lead to a precise description of $B^\times(P)$, when $P$ is a 2-group (actually also for arbitrary $p$-groups, but the case $p$ odd is known to be rather trivial). The main ingredient consists to show that $B^\times$ is a {\em rational} biset functor, and this is done by showing that the functor $B^\times$ (restricted to $p$-groups) is a subfunctor of the functor $\F_2R_\Q^*$. This leads to a description of $B^\times(P)$ in terms of a {\em genetic basis} of $P$, or equivalently, in terms of rational irreducible representations of $P$.\par
The paper is organized as follows~: in Section~2, I recall the main definitions and notation on biset functors. Section~3 deals with genetic subgroups and rational biset functors. Section~4 gives a natural exposition of the biset functor structure of $B^\times$. In Section~5, I state results about faithful elements in $B^\times(P)$ for some specific $p$-groups $P$. In Section~6, I introduce a natural transformation of biset functors from $B^\times$ to $\F_2B^*$. This transformation is injective, and in Section~7, I show that the image of its restriction to $p$-groups is contained in the subfunctor $\F_2R_\Q^*$ of $\F_2B^*$. This is the key result, leading in Section~8 to a description of the lattice of subfunctors of the restriction of~$B^\times$ to $p$-groups~: it is always a uniserial $p$-biset functor (even simple if $p$ is odd). This also provides an answer to the question, raised by Yal\c cin~(\cite{yalcin}), of the surjectivity of the exponential map $B(P)\to B^\times(P)$ for a 2-group~$P$.
\section{Biset functors}
\begin{enonce}{Notation and Definition} Denote by $\mathcal{C}$ the following category~:
\begin{itemize}
\item The objects of $\mathcal{C}$ are the finite groups.
\item If $G$ and $H$ are finite $p$-groups, then $\Hom_\mathcal{C}(G,H)=B(H\times G\op)$ is the Burnside group of finite $(H,G)$-bisets. An element of this group is called a {\em virtual} $(H,G)$-biset.
\item The composition of morphisms is $\Z$-bilinear, and if $G$, $H$, $K$ are finite groups, if $U$ is a finite $(H,G)$-biset, and $V$ is a finite $(K,H)$-biset, then the composition of (the isomorphism classes of) $V$ and $U$ is the (isomorphism class) of $V\times_HU$. The identity morphism $\Id_G$ of the group $G$ is the class of the set $G$, with left and right action by multiplication.
\end{itemize}
If $p$ is a prime number, denote by $\mathcal{C}_p$ the full subcategory of $\mathcal{C}$ whose objects are finite $p$-groups.\par
Let $\mathcal{F}$ denote the category of additive functors from $\mathcal{C}$ to the category $\gmod{\Z}$ of abelian groups. An object of $\mathcal{F}$ is called a {\em biset functor}. Similarly, denote by $\mathcal{F}_p$ the category of additive functors from $\mathcal{C}_p$ to $\gmod{\Z}$. An object of $\mathcal{F}_p$ will be called a {\em $p$-biset functor}.
\end{enonce}
If $F$ is an object of $\mathcal{F}$, if $G$ and $H$ are finite groups, and if $\varphi\in\hc{G}{H}$, then the image of $w\in F(G)$ by the map $F(\varphi)$ will generally be denoted by $\varphi(w)$. The composition $\psi\circ\varphi$ of morphisms $\varphi\in\hc{G}{H}$ and $\psi\in\hc{H}{K}$ will also be denoted by $\psi\times_H\varphi$.
\begin{enonce}{Notation} The {\em Burnside} biset functor (defined e.g. as the Yoneda functor $\hc{\un}{-}$), will be denoted by $B$. The functor of rational representations (see Section 1 of \cite{fonctrq}) will be denoted by $R_\Q$. The restriction of $B$ and $R_\Q$ to $\mathcal{C}_p$ will also be denoted by $B$ and $R_\Q$.
\end{enonce}
\subsection{Examples~:}\label{indresinfdefiso} Recall that this formalism of bisets gives a single framework for the usual operations of induction, restriction, inflation, deflation, and transport by isomorphism via the following correspondences~:
\begin{itemize}
\item If $H$ is a subgroup of $G$, then let $\Ind_H^G\in\hc{H}{G}$ denote the set~$G$, with left action of $G$ and right action of $H$ by multiplication.
\item If $H$ is a subgroup of $G$, then let $\Res_H^G\in\hc{G}{H}$ denote the set~$G$, with left action of $H$ and right action of $G$ by multiplication.
\item If $N\normal G$, and $H=G/N$, then let $\Inf_H^G\in\hc{H}{G}$ denote the set~$H$, with left action of $G$ by projection and multiplication, and right action of $H$ by multiplication.
\item If $N\normal G$, and $H=G/N$, then let $\Def_H^G\in\hc{G}{H}$ denote the
 set~$H$, with left action of $H$ by multiplication, and right action of $G$ by  projection and multiplication.
\item If $\varphi: G\to H$ is a group isomorphism, then let $\Iso_G^H=\Iso_G^H(\varphi)\in \hc{G}{H}$ denote the set~$H$, with left action of $H$ by multiplication, and right action of $G$ by taking image by $\varphi$, and then multiplying in~$H$.
\end{itemize}
\begin{enonce}{Definition} A {\em section} of the group $G$ is a pair $(T,S)$ of subgroups of $G$ such that $S\normal T$.
\end{enonce}
\begin{enonce}{Notation} If $(T,S)$ is a section of $G$, set
$$\Indinf_{T/S}^G=\Ind_T^G\Inf_{T/S}^T\ressort{1cm}\hbox{and}\ressort{1cm}\Defres_{T/S}^G=\Def_{T/S}^T\Res_T^G\mpoint$$
\end{enonce}
Then $\Indinf_{T/S}^G\cong G/S$ as $(G,T/S)$-biset, and $\Defres_{T/S}^G\cong S\dom G$ as $(T/S,G)$-biset.
\begin{enonce}{Notation} \label{Tu}Let $G$ and $H$ be groups, let $U$ be an $(H,G)$-biset, and let $u\in U$. If $T$ is a subgroup of $H$, set
$$T^u=\{g\in G\mid \exists t\in T,\; tu=ug\}\mpoint$$
This is a subgroup of $G$. Similarly, if $S$ is a subgroup of $G$, set
$${^u}S=\{h\in H\mid\exists s\in S,\;us=hu\}\mpoint$$
This is a subgroup of $H$.
\end{enonce}
\begin{enonce}{Lemma}\label{UmodG} Let $G$ and $H$ be groups, let $U$ be an $(H,G)$-biset, and let $S$ be a subgroup of $G$. Then there is an isomorphism of $H$-sets
$$U/G=\bigsqcup_{u\in [H\dom U/S]}H/{^uS}\mvirg$$
where $[H\dom U/S]$ is a set of representatives of $(H,S)$-orbits on $U$.
\end{enonce}
\pf Indeed $H\dom U/S$ is the set of orbits of $H$ on $U/S$, and $^uS$ is the stabilizer of~$uS$ in~$H$.
\subsection{Opposite bisets~:} If $G$ and $H$ are finite groups, and if $U$ is a finite $(H,G)$-biset, then let $U\op$ denote the opposite biset~: as a set, it is equal to~$U$, and it is a $(G,H)$-biset for the following action
$$\forall h\in H,\forall u\in U,\forall g\in G,\;g.u.h\;({\rm in}\;U\op)=h^{-1}ug^{-1}\;({\rm in}\;U)\mpoint$$
This definition can be extended by linearity, to give an isomorphism
$$\varphi\mapsto\varphi\op: \hc{G}{H}\to\hc{H}{G}\mpoint$$
It is easy to check that $(\varphi\circ\psi)\op=\psi\op\circ\varphi\op$, with obvious notation, and the functor
$$\left\{ \begin{array}{l}G\mapsto G\\\varphi\mapsto\varphi\op\end{array}\right.$$
is an equivalence of categories from $\mathcal{C}$ to the dual category, which restricts to an equivalence of $\mathcal{C}_p$ to its dual category.\par
\begin{rem}{Example} if $G$ is a finite group, and $(T,S)$ is a section of $G$, then
$$(\Indinf_{T/S}^G)\op\cong\Defres_{T/S}^G$$
as $(T/S,G)$-bisets.
\end{rem}
\begin{enonce}{Definition and Notation} If $F$ is a biset functor, the {\em dual} biset functor~$F^*$ is defined by
$$F^*(G)=\Hom_\Z(F(G),\Z)\mvirg$$
for a finite group $G$, and by
$$F^*(\varphi)(\alpha)=\alpha\circ F(\varphi\op)\mvirg$$
for any $\alpha\in F^*(G)$, any finite group $H$, and any $\varphi\in\hc{G}{H}$.
\end{enonce}

\subsection{Some idempotents in $\ec{G}$~:} Let $G$ be a finite group, and let $N\normal G$. Then it is clear from the definitions that
$$\Def_{G/N}^G\circ \Inf_{G/N}^G=(G/N)\times_G(G/N)=\Id_{G/N}\mpoint$$
It follows that the composition $e_N^G=\Inf_{G/N}^G\circ\Def_{G/N}^G$ is an idempotent in $\ec{G}$. Moreover, if $M$ and $N$ are normal subgroups of $G$, then $e_N^G\circ e_M^G=e_{NM}^G$. Moreover $e_1^G=\Id_G$. 
\begin{enonce}{Lemma}{\rm (\cite{bisetsections} Lemma 2.5)} \label{decomposition} If $N\normal G$, define $f_N^G\in\ec{G}$ by
$$f_N^G=\sumb{M\normal G}{N\subseteq M}\mu_{\normal G}(N,M)e_M^G\mvirg$$
where $\mu_{\normal G}$ denotes the M\"obius function of the poset of normal subgroups of~$G$. Then the elements $f_N^G$, for $N\normal G$, are orthogonal idempotents of $\ec{G}$, and their sum is equal to $\Id_G$.
\end{enonce}
Moreover, it is easy to check from the definition that for $N\normal G$,
\begin{equation}\label{somme partielle}
f_N^G=\Inf_{G/N}^G\circ f_\un^{G/N}\circ \Def_{G/N}^G\mvirg
\end{equation}
and
$$e_N^G=\Inf_{G/N}^G\circ \Def_{G/N}^G=\sumb{M\normal G}{M\supseteq N}f_M^G\mpoint$$
\begin{enonce}{Lemma}\label{partial}
If $N$ is a non trivial normal subgroup of $G$, then 
$$f_\un^G\circ \Inf_{G/N}^G=0\;\;\;\hbox{ and }\;\;\;\Def_{G/N}^G\circ f_\un^G=0\mpoint$$
\end{enonce}
\pf Indeed by~\ref{somme partielle}
\begin{eqnarray*}
f_\un^G\circ \Inf_{G/N}^G&=&f_\un^G\circ \Inf_{G/N}^G\circ \Def_{G/N}^G\circ\Inf_{G/N}^G\\
&=&\sumb{M\normal N}{M\supseteq N}f_\un^Gf_M^G\Inf_{G/N}^G=0\mvirg
\end{eqnarray*}
since $M\neq \un$ when $M\supseteq N$. The other equality of the lemma follows by taking opposite bisets.\findemo
\begin{rem}{Remark}\label{f1p}
It was also shown in Section 2.7~of~\cite{bisetsections} that if $P$ is a $p$-group, then
$$f_\un^P=\sum_{N\subseteq \Omega_1Z(P)}\mu(\un,N)P/N\mvirg$$
where $\mu$ is the M\"obius function of the poset of subgroups of~$N$, and $\Omega_1Z(P)$ is the subgroup of the centre of $P$ consisting of elements of order at most $p$.
\end{rem}
\begin{enonce}{Notation and Definition} If $F$ is a a biset functor, and if $G$ is a finite group, then the idempotent $f_\un^G$ of $\ec{G}$ acts on $F(G)$. Its image
$${\partial}F(G)=f_\un^GF(G)$$
is a direct summand of $F(G)$ as $\Z$-module~: it will be called the set of {\em faithful} elements of $F(G)$.
\end{enonce}
The reason for this name is that any element $u\in F(G)$ which is inflated from a proper quotient of $G$ is such that $F(f_\un^G)u=0$. From Lemma~\ref{partial}, it is also clear that
$$\partial F(G)=\mathop{\bigcap}_{\un\neq N\normal G}\Ker\;\Def_{G/N}^G\mpoint$$
\section{Genetic subgroups and rational $p$-biset functors}
The following definitions are essentially taken from Section~2 of~\cite{dadegroup}~:
\begin{enonce}{Definition and Notation} Let $P$ be a finite $p$-group. If $Q$ is a subgroup of $P$, denote by $Z_P(Q)$ the subgroup of $P$ defined by
$$Z_P(Q)/Q=Z(N_P(Q)/Q)\mpoint$$
A subgroup $Q$ of $P$ is called {\em genetic} if it satisfies the following two conditions~:
\begin{enumerate}
\item The group $N_P(Q)/Q$ has normal $p$-rank 1.
\item If $x\in P$, then $Q^x\cap Z_P(Q)\subseteq Q$ if and only if $Q^x=Q$.
\end{enumerate}
Two genetic subgroups $Q$ and $R$ are said to be {\em linked modulo $P$} (notation $Q\estliemod{P}R$), if there exist elements $x$ and $y$ in $P$ such that $Q^x\cap Z_P(R)\subseteq R$ and $R^y\cap Z_P(Q)\subseteq Q$.\par
This relation is an equivalence relation on the set of genetic subgroups of~$P$. The set of equivalence classes is in one to one correspondence with the set of isomorphism classes of rational irreducible representations of $P$. A {\em genetic basis} of $P$ is a set of representatives of these equivalences classes.
\end{enonce} 
If $V$ is an irreducible representation of $P$, then the {\em type} of $V$ is the isomorphism class of the group $N_P(Q)/Q$, where $Q$ is a genetic subgroup of $P$ in the equivalence class corresponding to $V$ by the above bijection.
\begin{rem}{Remark} The definition of the relation $\estliemod{P}$ given here is different from Definition~2.9 of~\cite{dadegroup}, but it is equivalent to it, by Lemma~4.5 of~\cite{bisetsections}.
\end{rem} 
The following is Theorem~3.2 of~\cite{bisetsections}, in a slightly different form~:
\begin{enonce}{Theorem} Let $P$ be a finite $p$-group, and $\mathcal{G}$ be a genetic basis of $P$. Let $F$ be a $p$-biset functor. Then the map
$$\mathcal{I}_\mathcal{G}=\oplus_{Q\in\mathcal{G}}\Indinf_{N_P(Q)/Q}^P:\oplus_{Q\in\mathcal{G}}\partial F\big(N_P(Q)/Q\big)\to F(P)$$
is split injective.
\end{enonce}
\begin{rem}{Remark} There are two differences with the initial statement of Theorem~3.2 of~\cite{bisetsections}~: here I use genetic {\em subgroups} instead of genetic {\em sections}, because these two notions are equivalent by Proposition~4.4 of~\cite{bisetsections}. Also the definition of the map $\mathcal{I}_\mathcal{G}$ is apparently different~: with the notation of~\cite{bisetsections}, the map $\mathcal{I}_\mathcal{G}$ is the sum of the maps $F(a_Q)$, where $a_Q$ is the trivial $(P,P/P)$-biset if $Q=P$, and $a_Q$ is the virtual $(P,N_P(Q)/Q)$-biset $P/Q-P/\hat{Q}$ if $Q\neq P$, where $\hat{Q}$ is the unique subgroup of $Z_P(Q)$ containing $Q$, and such that $|\hat{Q}:Q|=p$. But it is easy to see that the restriction of the map $F(P/\hat{Q})$ to $\partial F(N_P(Q)/Q)$ is actually~0. Moreover, the map $F(a_Q)$ is equal to $\Indinf_{N_P(Q)/Q}^P$. So in fact, the above map $\mathcal{I}_\mathcal{G}$ is the same as the one defined in Theorem~3.2 of~\cite{bisetsections}.
\end{rem}
\begin{enonce}{Definition} A $p$-biset functor $F$ is called {\em rational} if for any finite $p$-group $P$ and any genetic basis $\mathcal{G}$ of $P$, the map $\mathcal{I}_\mathcal{G}$ is an isomorphism.
\end{enonce}
It was shown in Proposition~7.4 of~\cite{bisetsections} that subfunctors, quotient functors, and dual functors of rational $p$-biset functors are rational.
\section{The functor of units of the Burnside ring}
\begin{enonce}{Notation} If $G$ is a finite group, let $B^\times(G)$ denote the group of units of the Burnside ring $B(G)$.
\end{enonce}
If $G$ and $H$ are finite groups, if $U$ is a finite $(H,G)$-biset, recall that $U\op$ denotes the $(G,H)$-biset obtained from $U$ by reversing the actions. If $X$ is a finite $G$-set, then $T_U(X)=\Hom_G(U\op,X)$ is a finite $H$-set. The correspondence $X\mapsto T_U(X)$ can be extended to a correspondence $T_U:B(G)\to B(H)$, which is multiplicative (i.e. $T_U(ab)=T_U(a)T_U(b)$ for any $a,b\in B(G)$), and preserves identity elements (i.e. $T_U(G/G)=H/H$). This extension to $B(G)$ can be built by different means, and the following is described in Section~4.1 of \cite{tensams}~: if $a$ is an element of $B(G)$, then there exists a finite $G$-poset $X$ such that $a$ is equal to the Lefschetz invariant $\Lambda_X$. Now $\Hom_G(U\op,X)$ has a natural structure of $H$-poset, and one can set $T_U(a)=\Lambda_{{\scriptstyle\rm Hom}_G(U\op,X)}$. It is an element of $B(H)$, which does not depend of the choice of the poset $X$ such that $a=\Lambda_X$, because with Notation~\ref{Tu} and Lemma~\ref{UmodG}, for any subgroup $T$ of $H$ the Euler-Poincar\'e characteristics $\chi\left(\Hom_G(U\op,X)^T\right)$ can be computed by
$$\chi\left(\Hom_G(U\op,X)^T\right)=\prod_{u\in T\dom U/G}\chi(X^{T^u})\mvirg$$
and the latter only depends on the element $\Lambda_X$ of $B(G)$. As a consequence, one has that
$$|T_U(a)^T|=\prod_{u\in T\dom U/G}|a^{T^u}|\mpoint$$
It follows in particular that $T_U\big(B^\times(G)\big)\subseteq B^\times(H)$. Moreover, it is easy to check that $T_U=T_{U'}$ if $U$ and $U'$ are isomorphic $(H,G)$-bisets, that $T_{U_1\sqcup U_2}(a)=T_{U_1}(a)T_{U_2}(a)$ for any $(H,G)$-bisets $U_1$ and $U_2$, and any $a\in B(G)$.\par
It follows that there is a well defined bilinear pairing
$$B(H\times G\op)\times B^\times(G)\to B^\times(H)\mvirg$$
extending the correspondence $(U,a)\mapsto T_U(a)$. If $f\in B(H\times G\op)$ (i.e. if $f$ is a virtual $(H,G)$-biset), the corresponding group homomorphism $B^\times(G)\to B^\times(H)$ will be denoted by $B^\times (f)$. \par
Now let $K$ be a third group, and $V$ be a finite $(K,H)$-set. If $X$ is a finite $G$-set, there is a canonical isomorphism of $K$-sets
$$\Hom_H\big(V\op,\Hom_G(U\op,X)\big)\cong \Hom_G\big((V\times_HU)\op,X\big)\mvirg$$
showing that $T_V\circ T_U=T_{V\times_HU}$.\par
It follows more generally that $B^\times(g)\circ B^\times(f)=B^\times(g\times_Hf)$ for any $g\in B(K\times H\op)$ and any $f\in B(H\times G\op)$. Finally this shows~:
\begin{enonce}{Proposition} The correspondence sending a finite group $G$ to $B^\times(G)$, and an homomorphism $f$ in $\mathcal{C}$ to $B^\times(f)$, is a biset functor.
\end{enonce}
\begin{rem}{Remark and Notation} The restriction and inflation maps for the functor $B^\times$ are the usual ones for the functor $B$. The deflation map $\Def_{G/N}^G$ corresponds to taking fixed points under $N$ (so it {\em does not coincide} with the usual deflation map for $B$, which consist in taking {\em orbits under $N$}). \par
Similarly, if $H$ is a subgroup of $G$, the induction map from $H$ to $G$ for the functor $B^\times$ is sometimes called {\em multiplicative induction}. I will call it {\em tensor induction}, and denote it by $\Ten_H^G$. If $(T,S)$ is a section of $G$, I will also set $\Teninf_{T/S}^P=\Ten_T^P\Inf_{T/S}^T$.
\end{rem}
\vspace{1cm}
\section{Faithful elements in $B^\times(G)$}
\begin{enonce}{Notation and definition} Let $G$ be a finite group. Denote by~$[s_G]$ a set of representatives of conjugacy classes of subgroups of $G$. Then the elements $G/L$, for $L\in [s_G]$, form a basis of $B(G)$ over $\Z$, called the {\em canonical basis} of $B(G)$.
\end{enonce}
The primitive idempotents of $\Q B(G)$ are also indexed by $[s_G]$~: if $H\in[s_G]$, the correspondent idempotent $e_H^G$ is equal to
$$e_H^G=\frac{1}{|N_G(H)|}\sum_{K\subseteq H}|K|\mu(K,H)G/K\mvirg$$
where $\mu(K,H)$ denotes the M\"obius function of the poset of subgroups of $G$, ordered by inclusion (see \cite{gluck}, \cite{yoshidaidemp}, or \cite{handbook}).\par
Recall that if $a\in B(G)$, then $a\cdot e_H^G=|a^H|e_H^G$ so that $a$ can be written as
$$ a=\sum_{H\in[s_G]}|a^H|e_H^G\mpoint$$
Now $a\in B^\times(G)$ if and only if $a\in B(G)$ and $|a^H|\in\{\pm 1\}$ for any $H\in [s_G]$, or equivalently if $a^2=G/G$. If now $P$ is a $p$-group, and if $p\neq 2$, since $|a^H|\equiv |a|\;(p)$ for any subgroup $|H|$ of $P$, it follows that $|a^H|=|a|$ for any $H$, thus $a=\pm P/P$. This shows the following well know
\begin{enonce}{Lemma}\label{impair} If $P$ is an odd order $p$-group, then $B^\times(P)=\{\pm P/P\}$.
\end{enonce}
\begin{rem}{Remark}
So in the sequel, when considering $p$-groups, the only really non-trivial case will occur for $p=2$. However, some statements will be given for arbitrary $p$-groups.
\end{rem}
\begin{enonce}{Notation} If $G$ is a finite group, denote by $F_G$ the set of subgroups~$H$ of~$G$ such that $H\cap Z(G)=\un$, and set $[F_G]=F_G\cap[s_G]$.
\end{enonce}
\begin{enonce}{Lemma}\label{groscentre} Let $G$ be a finite group. If $|Z(G)|>2$, then $\partial B^\times(G)$ is trivial.
\end{enonce}
\pf Indeed let $a\in\partial B^\times(G)$. Then $\Def_{G/N}^Ga$ is the identity element of $B^\times(G/N)$, for any non-trivial normal subgroup $N$ of $G$. Now suppose that $H$ is a subgroup of $G$ containing $N$. Then
$$|a^H|=|\Defres_{N_G(H)/H}^Ga|=|\Iso_{N_{G/N}(H/N)/(H/N)}^{N_G(H)/H}\Defres_{N_{G/N}(H/N)}^{G/N} \Def_{G/N}^Ga|=1\mpoint$$
In particular $|a^H|=1$ if $H\cap Z(G)\neq \un$. It follows that there exists a subset~$A$ of~$[F_G]$ such that
$$a=G/G-2\sum_{H\in A}e_H^G\mpoint$$
If $A\neq\emptyset$, i.e. if $a\neq G/G$, let $L$ be a maximal element of $A$. Then $L\neq G$, because $Z(G)\neq \un$. The coefficient of $G/L$ in the expression of $a$ in the canonical basis of $B(G)$ is equal to
$$-2\frac{|L|\mu(L,L)}{|N_G(L)|}=-2\frac{|L|}{|N_G(L)|}\mpoint$$
This is moreover an integer, since $a\in B^\times(G)$. It follows that $|N_G(L):L|$ is equal to 1 or 2. But since $L\cap Z(G)=\un$, the group $Z(G)$ embeds into the group $N_G(L)/L$. Hence $|N_G(L):L|\geq 3$, and this contradiction shows that $A=\emptyset$, thus $a=G/G$.\findemo
\begin{enonce}{Lemma} \label{critere} Let $P$ be a finite $2$-group, of order at least 4, and suppose that the maximal elements of $F_P$ have order 2. If $|P|\geq 2|F_P|$, then $\partial B^\times(P)$ is trivial.
\end{enonce}
\pf Let $a\in\partial B^\times(P)$. By the argument of the previous proof, there exists a subset $A$ of $[F_P]$ such that
$$a=P/P-2\sum_{H\in A}e_H^P\mpoint$$
The hypothesis implies that $\mu(\un,H)=-1$ for any non-trivial element $H$ of $[F_P]$. Now if $\un\in A$, the coefficient of $P/1$ in the expression of $a$ in the canonical basis of $B(P)$ is equal to 
$$-2\frac{1}{|P|}+2\sum_{H\in A-\{\un\}}\frac{1}{|N_P(H)|}=-2\frac{1}{|P|}+2\sum_{H\in \sur{A}-\{\un\}}\frac{1}{|P|}=\frac{-4+2|\sur{A}|}{|P|}\mvirg$$
where $\sur{A}$ is the set of subgroups of $P$ which are conjugate to some element of~$A$. This coefficient is an integer if $a\in B(P)$, so $|P|$ divides $2|\sur{A}|-4$. But~$|\sur{A}|$ is always odd, since the trivial subgroup is the only normal subgroup of $P$ which is in $\sur{A}$ in this case. Thus $2|\sur{A}|-4$ is congruent to 2 modulo 4, and cannot be divisible by $|P|$, since $|P|\geq 4$. \par
So $\un\notin A$, and the coefficient of $P/1$ in the expression of $a$ is equal to
$$2\sum_{H\in A}\frac{1}{|N_P(H)|}=\frac{2|\sur{A}|}{|P|}\mpoint$$
Now this is an integer, so $2|\sur{A}|$ is congruent to $0$ or $1$ modulo the order of~$P$, which is even since $|P|\geq 2|F_P|\geq 2$. Thus $\un\notin A$, and $2|\sur{A}|$ is a multiple of~$|P|$. But $2|\sur{A}|< 2|F_P|$ since $\un\notin A$. So if $2|F_P|\leq|P|$, it follows that
$\sur{A}$ is empty, and $A$ is empty. Hence $a=P/P$, as was to be shown.\findemo
\begin{enonce}{Corollary} \label{pasfidele}Let $P$ be a finite 2-group. Then the group $\partial B^\times(P)$ is trivial in each of the following cases~:
\begin{enumerate}
\item $P$ is abelian of order at least 3.
\item $P$ is generalized quaternion or semi-dihedral.
\end{enumerate}
\end{enonce}
\begin{rem}{Remark} Case 1 follows easily from Matsuda's Theorem (\cite{matsuda}). Case~2 follows from Lemma~4.6 of Yal\c cin (\cite{yalcin}). 
\end{rem}
\pf Case 1 follows from Lemma~\ref{groscentre}. In Case 2, if $P$ is generalized quaternion, then $F_P=\{\un\}$, thus $|P|\geq 2|F_P|$. And if $P$ is semidihedral, then there is a unique conjugacy class of non-trivial subgroups $H$ of~$P$ such that $H\cap Z(P)=\un$. Such a group has order 2, and $N_P(H)=HZ(P)$ has order~4. Thus $|F_P|=1+\frac{|P|}{4}$, and $|P|\geq ²2|F_P|$ also in this case.\findemo
\begin{enonce}{Corollary} {\rm [Yal\c cin \cite{yalcin} Lemma~4.6 and Lemma~5.2]} \label{upsilon} Let $P$ be a $p$-group of normal $p$-rank 1. Then $\partial B^\times(P)$ is trivial, except if $P$ is 
\begin{itemize}
\item the trivial group, and $\partial B^\times(P)$ is the group of order 2 generated by $\upsilon_P=-P/P$.
\item cyclic of order 2, and $\partial B^\times(P)$ is the group of order 2 generated by 
$$\upsilon_P=P/P-P/\un\mpoint$$
\item dihedral of order at least 16, and then $\partial B^\times(P)$ is the group of order 2 generated by the element 
$$\upsilon_P=P/P+P/1-P/I-P/J\mvirg$$
where $I$ and $J$ are non-central subgroups of order 2 of $P$, not conjugate in $P$.
\end{itemize}
\end{enonce}
\pf Lemma~\ref{impair} and Lemma~\ref{groscentre} show that $\partial B^\times(P)$ is trivial, when $P$ has normal $p$-rank 1, and $P$ is not trivial, cyclic of order 2, or dihedral~: indeed then, the group $P$ is cyclic of order at least 3, or generalized quaternion, or semi-dihedral.\par
Now if $P$ is trivial, then obviously $B(P)=\Z$, so $B^\times(P)=\partial B^\times(P)=\{\pm P/P\}$. If $P$ has order 2, then clearly $B^\times(P)$ consists of $\pm P/P$ and $\pm(P/P-P/\un)$, and $\partial B^\times(P)=\{P/P,P/P-P/\un\}$. Finally, if $P$ is dihedral, the set $F_P$ consists of the trivial group, and of two conjugacy classes of subgroups $H$ of order~2 of~$P$, and $N_P(H)=HZ$ for each of these, where~$Z$ is the centre of $P$. Thus 
$$|F_P|=1+2\frac{|P|}{4}=1+\frac{|P|}{2}\mpoint$$
Now with the notation of the proof of Lemma~\ref{critere}, one has that $2|\sur{A}|\equiv 0\;(|P|)$, and $2|\sur{A}|<|F_P|=2+|P|$. So either $A=\emptyset$, and in this case $a=P/P$, or $2|\sur{A}|=|P|$, which means that $\sur{A}$ is the whole set of non-trivial elements of~$F_P$. In this case
$$a=P/P-2(e_I^P+e_J^P)\mvirg$$
where $I$ and $J$ are non-central subgroups of order 2 of $P$, not conjugate in~$P$. It is then easy to check that
$$a=P/P+P/\un-(P/I+P/J)\mvirg$$
so $a$ is indeed in $B(P)$, hence in $B^\times(P)$. Moreover $\Def_{P/Z}^Pa$ is the identity element of $B^\times(P/Z)$, so $a=f_\un^Pa$, and $a\in\partial B^\times(P)$. This completes the proof.\findemo

\section{A morphism of biset functors}
If $k$ is any commutative ring, there is an obvious isomorphism of biset functor from $kB^*=k\otimes_\Z B^*$ to $\Hom(B,k)$, which is defined for a group~$G$ by sending the element $\alpha=\sum_i\alpha_i\otimes \psi_i$, where $\alpha_i\in k$ and $\psi_i\in B^*(G)$, to the linear form $\tilde{\alpha}:B(G)\to k$ defined by $\tilde{\alpha}(G/H)=\sum_i \psi_i(G/H)\alpha_i$.
\begin{enonce}{Notation} Let $\{\pm 1\}=\Z^\times$ be the group of units of the ring $\Z$. The unique group isomorphism from $\{\pm 1\}$ to $\Z/2\Z$ will be denoted by $u\mapsto u_+$.
\end{enonce}
If $G$ is a finite group, and if $a\in B^\times(G)$, then recall that for each subgroup~$S$ of $G$, the integer $|a^S|$ is equal to $\pm 1$. Define a map $\epsilon_G : B^\times(G)\to \F_2B^*(G)$ by setting $\epsilon_G(a)(G/S)=|a^S|_+$, for any $a\in B^\times(G)$ and any subgroup $S$ of $G$.
\begin{enonce}{Proposition} The maps $\epsilon_G$ define a injective morphism of biset functors $$\epsilon:B^\times \to \F_2B^*\mpoint$$
\end{enonce}
\pf The injectivity of the map $\epsilon_G$ is obvious. Now let $G$ and $H$ be finite groups, and let $U$ be a finite $(H,G)$-biset. Also denote by $U$ the corresponding element of $B(H\times G\op)$. If $a\in B^\times(G)$, and if $T$ is a subgroup of $H$, then
$$|B^\times(U)(a)^T|=\prod_{u\in T\dom U/G}|a^{T^u}|\mpoint$$
Thus
\begin{eqnarray*}
\epsilon_H\left(B^\times(U)(a)\right)(H/T)&=&\left(\prod_{u\in T\dom U/G}|a^{T^u}|\right)_+\\
&=&\sum_{u\in T\dom U/G}|a^{T^u}|_+\\
&=&\sum_{u\in T\dom U/G}\epsilon_G(a)(G/T^u)\\
&=&\epsilon_G(a)(U\op/T)\\
&=&\epsilon_G(a)(U\op\times_HH/T)\\
&=&\F_2B^*(U)\big(\epsilon_G(a)\big)(H/T)
\end{eqnarray*}
thus $\epsilon_H\circ B^\times(U)=\F_2B^*(U)\circ\epsilon_G$. Since both sides are additive with respect to $U$, the same equality holds when $U$ is an arbitrary element of $B(H\times G\op)$, completing the proof.\findemo
\section{Restriction to $p$-groups}
The additional result that holds for finite $p$-groups (and not for arbitrary finite groups) is the Ritter-Segal theorem, which says that the natural transformation $B\to R_\Q$ of biset functors for $p$-groups, is surjective. By duality, it follows that the natural transformation $i:kR_\Q^*\to kB^*$ is injective, for any commutative ring $k$. The following gives a characterization of the image $i(kR_\Q^*)$ inside $kB^*$~:
\begin{enonce}{Proposition} \label{caract}Let $p$ be a prime number, let $P$ be a $p$-group, let $k$ be a commutative ring.Then the element $\varphi\in kB^*(P)$ lies in $i\big(kR_\Q^*(P)\big)$ if and only if the element $\Defres_{T/S}^P\varphi$ lies in $i\big(kR_\Q^*(T/S)\big)$, for any section $T/S$ of~$P$ which is 
\begin{itemize}
\item elementary abelian of rank~2, or non-abelian of order $p^3$ and exponent~$p$, if $p\neq 2$.
\item elementary abelian of rank~2, or dihedral of order at least 8, if $p=2$.
\end{itemize}
\end{enonce} 
\pf Since the image of $kR_\Q^*$ is a subfunctor of $kB^*$, if $\varphi\in i\big(kR_\Q^*(P)\big)$, then $\Defres_{T/S}^P\varphi\in i\big(kR_\Q^*(T/S)\big)$, for any section $(T,S)$ of $P$. \par
Conversely, consider the exact sequence of biset functors over $p$-groups
$$0\to K\to B\to R_Q\to 0\mpoint$$
Every evaluation of this sequence at a particular $p$-group is a split exact sequence of (free) abelian groups. Hence by duality, for any ring $k$, there is an exact sequence
$$0\to kR_\Q^*\to kB^*\to kK^*\to 0\mpoint$$
With the identification $kB^*\cong\Hom_\Z(B,k)$, this means that if $P$ is a $p$-group, the element $\varphi\in RB^*(P)$ lies in $i\big(kR_\Q^*(P)\big)$ if and only if $\varphi\big(K(P)\big)=0$. Now by Corollary~6.16 of \cite{dadegroup}, the group $K(P)$ is the set of linear combinations of elements of the form $\Indinf_{T/S}^P\theta(\kappa)$, where $T/S$ is a section of $P$, and $\theta$ is a group isomorphism from one of the group listed in the proposition to $T/S$, and $\kappa$ is a specific element of $K(T/S)$ in each case. The proposition follows, because
$$\varphi\big(\Indinf_{T/S}^P\theta(\kappa)\big)=(\Defres_{T/S}^P\varphi)\big(\theta(\kappa)\big)\mvirg$$
and this is zero if $\Defres_{T/S}^P\varphi$ lies in $i\big(kR_\Q^*(T/S)\big)$.\findemo
\begin{enonce}{Theorem} Let $p$ be a prime number, and $P$ be a finite $p$-group. The image of the map $\epsilon_P$ is contained in $i\big(\F_2R_\Q^*(P)\big)$.
\end{enonce}
\pf Let $a\in B^\times(P)$, and let $T/S$ be any section of $P$. Since
$$\Defres_{T/S}^Pi_P(a)=i_{T/S}\Defres_{T/S}^Pa\mvirg$$
by Proposition~\ref{caract}, it is enough to check that the image of $\epsilon_P$ is contained in~$i\big(\F_2R_\Q^*(P)\big)$, when $P$ is elementary abelian of rank~2 or non-abelian of order $p^3$ and exponent~$p$ if $p$ is odd, or when $P$ is elementary abelian of rank~2 or dihedral if $p=2$.\par
Now if $N$ is a normal subgroup of $P$, one has that 
$$f_N^Pi_P(a)=\Inf_{P/N}^P\left(i_{P/N}(f_{\un}^{P/N}\Def_{P/N}^Pa)\right)\mpoint$$
Thus by induction on the order of $P$, one can suppose $a\in \partial B^\times(P)$. But if $P$ is elementary abelian of rank~2, or if $P$ has odd order, then $\partial B^\times(P)$ is trivial, by Lemma~\ref{impair} and Corollary~\ref{pasfidele}. Hence there is nothing more to prove if $p$ is odd. And for $p=2$, the only case left is when $P$ is dihedral. In that case by Corollary~\ref{upsilon}, the group $\partial B^\times(P)$ has order 2, generated by the element
$$\upsilon_P=\sum_{H\in[s_P]-\{I,J\}}e_H^P-(e_I^P+e_J^P)\mvirg$$
where $[s_P]$ is a set of representatives of conjugacy classes of subgroups of $P$, and where $I$ and $J$ are the elements of $[s_P]$ which have order 2, and are non central in $P$. Moreover the element $\theta(\kappa)$ mentioned above is equal to
$$(P/I'-P/I'Z)-(P/J'-P/J'Z)\mvirg$$
where $Z$ is the centre of $P$, and $I'$ and $J'$ are non-central subgroups of order~2 of $P$, not conjugate in $P$. Hence up to sign $\theta(\kappa)$ is equal to 
$$\delta_P=(P/I-P/IZ)-(P/J-P/JZ)\mpoint$$ 
Since $\epsilon_P(\upsilon_P)(P/H)$ is equal to zero, except if $H$ is conjugate to $I$ or $J$, and then $\epsilon_P(\upsilon_P)(P/H)=1$, it follows that $\epsilon_P(\upsilon_P)(\delta_P)=1-1=0$, as was to be shown. This completes the proof.\findemo
\begin{enonce}{Corollary} The $p$-biset functor $B^\times$ is rational.
\end{enonce}
\pf Indeed, it is isomorphic to a subfunctor of $\F_2R_\Q^*\cong\Hom_\Z(R_\Q,\F_2)$, which is rational by Proposition~7.4 of~\cite{bisetsections}.\findemo
\begin{enonce}{Theorem} \label{units}Let $P$ be a $p$-group. Then $B^\times(P)$ is an elementary abelian 2-group of rank equal to the number of isomorphism classes of rational irreducible representations of~$P$ whose type is trivial, cyclic of order~2, or dihedral. More precisely~:
\begin{enumerate}
\item If $p\neq 2$, then $B^\times(P)=\{\pm 1\}$.
\item If $p=2$, then let $\mathcal{G}$ be a genetic basis of $P$, and let $\mathcal{H}$ be the subset of $\mathcal{G}$ consisting of elements $Q$ such that $N_P(Q)/Q$ is trivial, cyclic of order~2, or dihedral. If $Q\in\mathcal{H}$, then $\partial B^\times\big(N_P(Q)/Q\big)$ has order~2, generated by $\upsilon_{N_P(Q)/Q}$. Then the set
$$\{\Teninf_{N_P(Q)/Q}^P\upsilon_{N_P(Q)/Q}\mid Q\in \mathcal{H}\}$$
is an $\F_2$-basis of $B^\times(P)$.
\end{enumerate}
\end{enonce}
\pf This follows from the definition of a rational biset functor, and from Corollary~\ref{upsilon}. \findemo
\begin{rem}{Remark} If $P$ is abelian, then there is a unique genetic basis of $P$, consisting of subgroups $Q$ such that $P/Q$ is cyclic. So in that case, the rank of $B^\times(P)$ is equal 1 plus the number of subgroups of index 2 in $P$~: this gives a new proof of Matsuda's Theorem~(\cite{matsuda}).
\end{rem}
\section{The functorial structure of $B^\times$ for $p$-groups}
In this section, I will describe the lattice of subfunctors of the $p$-biset functor $B^\times$. 
\subsection{The case $p\neq 2$.} If $p\neq 2$, there is not much to say, since $B^\times(P)\cong \F_2$ for any $p$-group $P$. In this case, the functor $B^\times$ is the constant functor~$\Gamma_{\F_2}$ introduced in Corollary~8.4 of \cite{both}. It is also isomorphic to the simple functor~$S_{\un,\F_2}$. In this case, the results of~\cite{bisetsections} and~\cite{dadegroup} lead to the following remarkable version of Theorem~11.2 of~\cite{both}: 
\begin{enonce}{Proposition}\label{suiteexacte} If $p\neq 2$, the inclusion $B^\times\to \F_2R_\Q^*$ leads to a short exact sequence of $p$-biset functors
$$0\to B^\times\to \F_2R_\Q^*\to D_{tors}\to 0\mvirg$$
where $D_{tors}$ is the torsion part of the Dade $p$-biset functor.
\end{enonce}
\subsection{The case $p=2$.} There is a bilinear pairing 
$$\scal{\phantom{a}}{\phantom{b}}:\F_2R_\Q^*\times \F_2R_\Q\to\F_2\mpoint$$
This means that for each 2-group $P$, there is a bilinear form 
$$\scal{\phantom{a}}{\phantom{b}}_P: \F_2R_\Q^*(P)\times \F_2R_\Q(P)\to\F_2\mvirg$$
with the property that for any 2-group $Q$, for any $f\in\hcp{P}{Q}$, for any $a\in \F_2R_\Q^*(P)$ and any $b\in \F_2R_\Q(Q)$, one has that
$$\scal{\F_2R_\Q^*(f)(a)}{b}_Q=\scal{a}{\F_2R_\Q(f\op)(b)}_P\mpoint$$
Moreover this pairing is non-degenerate~: this means that for any 2-group~$P$, the pairing $\scal{\phantom{a}}{\phantom{b}}_P$ is non-degenerate. In particular, each subfunctor $F$ of $\F_2R_\Q^*$ is isomorphic to $\F_2R_\Q/F^\perp$, where $F^\perp$ is the orthogonal of $F$ for the pairing $\scal{\phantom{a}}{\phantom{b}}$.\par
In particular, the lattice of subfunctors of $\F_2R_\Q^*$ is isomorphic to the opposite lattice of subfunctors of $\F_2R_\Q$. Now since $B^\times$ is isomorphic to a subfunctor of $\F_2R_\Q$, its lattice of subfunctors is isomorphic to the opposite lattice of subfunctors of $\F_2R_\Q$ {\em containing} $B^\sharp=(B^\times)^\perp$. By Theorem~4.4 of \cite{fonctrq}, any subfunctor $L$ of $\F_2R_\Q$ is equal to the sum of subfunctors $H_Q$ it contains, where $Q$ is a 2-group of normal 2-rank 1, and $H_Q$ is the subfunctor of $\F_2R_\Q$ generated by the image $\sur{\Phi}_Q$ of the unique (up to isomorphism) irreducible rational faithful $\Q Q$-module $\Phi_Q$ in $\F_2R_\Q$. \par
In particular $B^\sharp$ is the sum of the subfunctors $H_Q$, where $Q$ is a 2-group of normal 2-rank 1 such that $\sur{\Phi}_Q\in B^\sharp(Q)$. This means that
$\scal{a}{\sur{\Phi}_Q}_Q=0$, for any $a\in B^\times(Q)$. Now $\Phi_Q=f_\un\Phi_Q$ since $\Phi_Q$ is faithful, so
$$\scal{a}{\sur{\Phi}_Q}_Q=\scal{a}{f_\un^Q\sur{\Phi}_Q}_Q=\scal{f_\un^Qa}{\sur{\Phi}_Q}_Q\mvirg$$
because $f_\un^Q=(f_\un^Q)\op$. Thus $\sur{\Phi}_Q\in B^\sharp(Q)$ if and only if $\sur{\Phi}_Q$ is orthogonal to $\partial B^\times(Q)$. Since $Q$ has normal $2$-rank 1, this is always the case by Corollary~\ref{upsilon}, except maybe if $Q$ is trivial, cyclic of order 2, or dihedral (of order at least 16). Now $H_\un=H_{C_2}=\F_2R_\Q$ by Theorem~5.6 of~\cite{fonctrq}. Since $B^\times$ is not the zero subfunctor of $\F_2R_\Q$, it follows that $H_Q\not\subseteq B^\sharp$, if $Q$ is trivial or cyclic of order 2. Now if $Q$ is dihedral, then $\Phi_Q$ is equal to $\Q Q/I-\Q Q/IZ$, where~$I$ is a non-central subgroup of order 2 of $Q$, and $Z$ is the centre of $Q$. Now
$$\epsilon_Q(\upsilon_Q)\big(i(\sur{\Phi}_Q)\big)=\epsilon_Q(\upsilon_Q)(Q/I-Q/IZ)=1-0=1\mvirg$$
It follows that $H_Q\not\subseteq B^\sharp$ if $Q$ is dihedral. Finally $B^\sharp$ is the sum of all subfunctors $H_Q$, when $Q$ is cyclic of order at least 4, or generalized quaternion, or semi-dihedral.\par
Recall from Theorem~6.2 of~\cite{fonctrq} that the poset of proper subfunctors of $\F_2R_\Q$ is isomorphic to the poset of closed subsets of the following graph~:
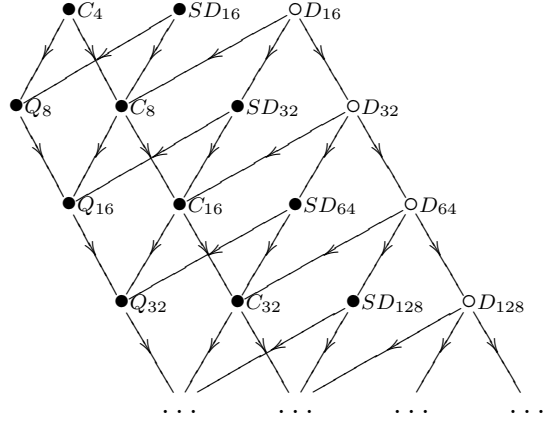
$$\xymatrix@C=10pt@R=30pt@M=0pt@H=0pt@W=0pt{
&\  \boite{C_4}\ & & \boite{SD_{16}}& & \boiteb{D_{16}}& & &\\
\boite{Q_8}\mar[ur]\mar[urrr]& & \boite{C_{8}}\mar[ul]\mar[ur]\mar[urrr]& & \boite{SD_{32}}\mar[ur]& & \boiteb{D_{32}}\mar[ul]& &\\
& \boite{Q_{16}}\mar[ul]\mar[ur]\mar[urrr]& & \boite{C_{16}}\mar[ul]\mar[ur]\mar[urrr]& & \boite{SD_{64}}\mar[ur]& & \boiteb{D_{64}}\mar[ul]& \\
& & \boite{Q_{32}}\mar[ul]\mar[ur]\mar[urrr]& & \boite{C_{32}}\mar[ul]\mar[ur]\mar[urrr]& & \boite{SD_{128}}\mar[ur]& & \boiteb{D_{128}}\mar[ul]\\
& & &{\rule{0pt}{8pt}\ldots\rule{0pt}{8pt}}\mar[ul]\mar[ur]\mar[urrr]& &{\rule{0pt}{8pt}\ldots\rule{0pt}{8pt}}\mar[ul]\mar[ur]\mar[urrr]& &{\rule{0pt}{8pt}\ldots\rule{0pt}{8pt}}\mar[ur]& &{\rule{0pt}{8pt}\ldots\rule{0pt}{8pt}}\mar[ul]\\
}$$
The vertices of this graph are the isomorphism classes of groups of normal 2-rank~1 and order at least~4, and there is an arrow from vertex $Q$ to vertex~$R$ if and only if $H_R\subseteq H_Q$. The vertices with a filled $\bullet$ are exactly labelled by the groups $Q$ for which $H_Q\subseteq B^\sharp$, and the vertices with a $\circ$ are labelled by dihedral groups.\par
By the above remarks, the lattice of subobjects of $B^\times$ is isomorphic to the opposite lattice of subfunctors of $\F_2R_\Q$ containing $B^\sharp$. Thus~:
\begin{enonce}{Theorem} The $p$-biset functor $B^\times$ is uniserial. It has an infinite strictly increasing series of proper subfunctors
$$0\subset L_0\subset L_1\cdots\subset L_n\subset\cdots$$
where $L_{0}$ is generated by the element $\upsilon_{\un}$, and $L_i$, for $i>0$, is generated by the element $\upsilon_{D_{2^{i+3}}}$ of $B^\times(D_{2^{i+3}})$. The functor $L_0$ is isomorphic to the simple functor $S_{\un,\F_2}$, and the quotient $L_{i}/L_{i-1}$, for $i\geq 1$, is isomorphic to the simple functor $S_{D_{2^{i+3}},\F_2}$.
\end{enonce}
\pf Indeed $L_0^\perp=B^\sharp+H_{D_{16}}$ is the unique maximal proper subfunctor of $\F_2R_\Q$. Thus $L_0$ is isomorphic to the unique simple quotient of $\F_2R_\Q$, which is $S_{\un,\F_2}$ by Proposition~5.1 of~\cite{fonctrq}. Similarly for $i\geq 1$, the simple quotient $L_i/L_{i-1}$ is isomorphic to the quotient
$$(B^\sharp+H_{D_{2^{i+3}}})/B^\sharp+H_{D_{2^{i+4}}})\mvirg$$
which is a quotient of 
$$(B^\sharp+H_{D_{2^{i+3}}})/B^\sharp\cong H_{D_{2^{i+3}}}/(B^\sharp\cap H_{D_{2^{i+3}}})\mpoint$$
But the only simple quotient of $H_{D_{2^{i+3}}}$ is $S_{D_{2^{i+3}},\F_2}$, by Proposition~5.1 of~\cite{fonctrq} again.\findemo
\begin{rem}{Remark} Let $P$ be a 2-group. By Theorem~5.12 of~\cite{fonctrq}, the $\F_2$-dimension of $S_{\un,\F_2}(P)$ is equal to the number of isomorphism classes of rational irreducible representations of~$P$ whose type is $\un$ or $C_2$, whereas the $\F_2$-dimension of $S_{D_{2^{i+3}},\F_2}(P)$ is the number of isomorphism classes of rational irreducible representations of~$P$ whose type is isomorphic to $D_{2^{i+3}}$. This gives a way to recover Theorem~\ref{units}~: the $\F_2$-dimension of $B^\times(P)$ is equal to the number of isomorphism classes of rational irreducible representations of~$P$ whose type is trivial, cyclic of order 2, or dihedral.
\end{rem}
\subsection{The surjectivity of the exponential map.} Let $G$ be a finite group. The exponential map ${\rm exp}_G: B(G)\to B^\times(G)$ is defined in Section~7 of Yal\c cin's paper~(\cite{yalcin}) by
$${\rm exp}_G(x)=(-1)\uparrow x\mvirg$$
where $-1=-\un/\un\in B^\times(\un)$, and where the exponentiation
$$(y,x)\in B^\times(G)\times B(G)\to B^\times(G)$$
is defined by extending the usual exponential map $(Y,X)\mapsto Y^X$, where $X$ and $Y$ are $G$-sets, and $Y^X$ is the set of maps from $X$ to $Y$, with $G$-action given by $(g\cdot f)(x)=gf(g^{-1}x)$.\par
Its possible to give another interpretation of this map~: indeed $B(G)$ is naturally isomorphic to $\hc{\un}{G}$, by considering any $G$-set as a $(G,\un)$-biset. It is clear that if $X$ is a finite $G$-set, and $Y$ is a finite set, then
$$T_X(Y)=Y^X\mpoint$$
This can be extended by linearity, to show that for any $x\in B(G)$
$$(-1)^x=B^\times(x)(-1)\mpoint$$
In particular the image $\Im({\rm exp}_G)$ of the exponential map ${\rm exp}_G$ is equal to $\hc{\un}{G}(-1)$. Denoting by $I$ the sub-biset functor of $B^\times$ generated by $-1\in B^\times(\un)$, it it now clear that 
$\Im({\rm exp}_G)=I(G)$ for any finite group $G$.\par
Now the restriction of the functor $I$ to the category $\mathcal{C}_2$ is equal to $L_0$, which is isomorphic to the simple functor $S_{\un,\F_2}$. Using Remark~5.13 of~\cite{fonctrq}, this shows finally the following~:
\begin{enonce}{Proposition} Let $P$ be a finite 2-group. Then~:
\begin{enumerate}
\item The $\F_2$-dimension of the image of the exponential map
$${\rm exp}_P: B(P)\to B^\times(P)$$
is equal to the number of isomorphism classes of absolutely irreducible rational representations of $P$.
\item The map ${\rm exp}_P$ is surjective if and only if the group $P$ has no irreducible rational representation of dihedral type, or equivalently, no genetic subgroup $Q$ such that $N_P(Q)/Q$ is dihedral.
\end{enumerate}
\end{enonce}
\begin{enonce}{Proposition} Let $p$ be a prime number. There is an exact sequence of $p$-biset functors~:
$$0\to B^\times\to \F_2R_\Q^*\to \F_2D^\Omega_{tors}\to 0\mvirg$$
where $D^\Omega_{tors}$ is the torsion part of the functor $D^\Omega$ of relative syzygies in the Dade group.
\end{enonce}
\pf In the case $p\neq 2$, this proposition is equivalent to Proposition~\ref{suiteexacte}, because $\F_2D^\Omega_{tors}=\F_2D_{tors}\cong D_{tors}$ in this case. And for $p=2$, the 2-functor $D^\Omega_{tors}$ is a quotient of the functor $R_\Q^*$, by Corollary~7.5 of~\cite{bisetsections}~: there is a surjective map $\pi: R_\Q^*\to D_{tors}^\Omega$, which is the restriction to $R_\Q^*$ of the surjection $\Theta: B^*\to D^\Omega$ introduced in Theorem~1.7 of~\cite{dadeburnside}. The $\F_2$-reduction of $\pi$ is a surjective map
$$\F_2\pi : \F_2R_\Q^*\to \F_2D^\Omega_{tors}\mpoint$$
To prove the proposition in this case, is is enough to show that the image of~$B^\times$ in $\F_2R_\Q^*$ is contained in the kernel of $\F_2\pi$, and that for any $2$-group~$P$, the $\F_2$-dimension of $\F_2R_\Q^*(P)$ is equal to the sum of the $\F_2$-dimensions of $B^\times(P)$ and $\F_2D^\Omega_{tors}(P)$~: but by Corollary~7.6 of~\cite{bisetsections}, there is a group isomorphism
$$D^\Omega_{tors}(P)\cong (\Z/4\Z)^{a_P}\oplus (\Z/2\Z)^{b_P}\mvirg$$
where $a_P$ is equal to the number of isomorphism classes of rational irreducible representations of~$P$ whose type is generalized quaternion, and $b_P$ equal to the number of isomorphism classes of rational irreducible representations of~$P$ whose type is cyclic of order at least 3, or semi-dihedral. Thus
$$\dim_{\F_2} \F_2D^\Omega_{tors}(P)=a_P+b_P\mpoint$$
Now since $\dim_{\F_2}B^\times(P)$ is equal to the number of isomorphism classes of rational irreducible representations of~$P$ whose type is cyclic of order at most 2, or dihedral, it follows that $\dim_{\F_2} \F_2D^\Omega_{tors}(P)+\dim_{\F_2}B^\times(P)$ is equal to the number of isomorphism classes of rational irreducible representations of $P$, i.e. to $\dim_{\F_2}\F_2R_\Q^*(P)$.\par
So the only thing to check to complete the proof, is that the image of~$B^\times$ in $\F_2R_\Q^*$ is contained in the kernel of $\F_2\pi$. Since $B^\times$, $\F_2R_\Q^*$ and $\F_2D^\Omega_{tors}$ are rational 2-biset functors, it suffices to check that if $P$ is a 2-group of normal 2-rank 1, and $a\in \partial B^\times(P)$, then the image of $a$ in $\partial\F_2R_\Q^*(P)$ lies in the kernel of $\F_2\pi$. There is nothing to do if $P$ is generalized quaternion, or semi-dihedral, or cyclic of order at least 3, for in this case $\partial B^\times(P)=0$ by Corollary~\ref{pasfidele}. Now if $P$ is cyclic of order at most 2, then $D^\Omega(P)=\zero$, and the result follows. And if $P$ is dihedral, then $D^\Omega(P)$ is torsion free by Theorem~10.3 of~\cite{cath}, so $D^\Omega_{tors}(P)=\zero$ again. \findemo

\end{document}